\documentclass{amsart}

\newtheorem{thm}{Theorem}

\newtheorem{lem}[thm]{Lemma}

\newcommand{\Pf}{\noindent{\bf Proof: }}

\begin{document}

\title{Bounds on the Crosscap Number of Torus Knots}       
\author{Thomas W.\ Mattman and Owen Sizemore}        
\address{Department of Mathematics and Statistics,
         California State University, Chico,
         Chico CA 95929-0525, USA}
\email{tmattman@csuchico.edu}
\address{Department of Mathematics, UC Berkeley,
	 970 Evans Hall \#3840,
	 Berkeley, CA 94720-3840, USA}
\email{owen\_s@uclink.berkeley.edu}
\thanks{The second author is an undergraduate student who was supervised 
by the first author during an REU held at CSU, Chico in
the summer of 2004 and funded by NSF REU Award 0354174.}

\begin{abstract}
For a torus knot $K$, we bound the crosscap number $c(K)$ in terms of
the genus $g(K)$ and crossing number $n(K)$: $c(K) \leq \lfloor (g(K)+9)/6 \rfloor$
and $c(K) \leq \lfloor (n(K) + 16)/12 \rfloor$. The $(6n-2,3)$ torus knots
show that these bounds are sharp.
\end{abstract}

\maketitle

\section{Introduction}
In 1978, Clark~\cite{C} defined the crosscap number $c(K)$ of the knot $K$ to be
the minimal genus of all non-orientable surfaces which span the knot and gave an
upper bound for this number in terms of $g(K)$, the genus: $c(K) \leq 2g(K) + 1$.
The obvious next question is if there is some way of bounding the genus in terms
of the crosscap number. Or, as Adams~\cite{A} asked, is there
some family of knots for which the difference $|g(K) - c(K)|$ increases without
bound? 

The torus knots are a natural target for this question since
Teragaito~\cite{T} has recently classified their crosscap
numbers. We soon noticed that the $(2n+1,2)$ torus knots
provide an answer to Adams's question. The genus of such a knot is $n$ while the
crosscap number is one. (Indeed, it's quite easy to see that these knots span
M\"obius bands. Take a strip of paper and give it $2n+1$ half twists before joining
the ends. The band's edge will be a $(2n+1,2)$ torus knot.) Thus,
for the $(2n+1,2)$ torus knots the difference  $g(K) - c(K)$ is $n-1$ and this
difference increases without bound as $n$ approaches infinity.

This example suggests that $2g(K)+1$ is a rather poor estimate for crosscap number
if we restrict attention to the class of torus knots. Indeed, we have the
following:
\begin{thm} 
 For a torus knot $K$, $c(K) \leq \lfloor (g(K)+9)/6 \rfloor$.
\label{propgen}
\end{thm}
\noindent%
Here, $\lfloor x \rfloor$ is the greatest integer less than or equal to $x$.

Since the crossing number of a torus knot is roughly twice the genus, we can
also improve the bound on crosscap number in terms of the
crossing number $n(K)$.  
\begin{thm} For a torus knot $K$, $c(K) \leq \lfloor (n(K) + 16)/12 \rfloor$.
\label{propcross}
\end{thm}
\noindent%
Compare this with Murakami and 
Yasuhara's~\cite{MY} general result that $c(K) \leq \lfloor n(K)/2
\rfloor$.

The $(6n-2,3)$ torus knots show that the inequalities in our two theorems are
sharp. These knots have genus $6n-3$ and crossing number $12n-4$. Thus,
both $(g(K)+9)/6$ and $(n(K) + 16)/12$ yield the crosscap number $n+1$ for these
knots.

Our arguments make use of Teragaito's~\cite{T} classification of the crosscap
number of torus knots. He shows that the crosscap number of a $(p,q)$ torus 
knot is given by summing certain coefficients in the continued fraction
expansion of $p/q$, $q/p$, $(pq+1)/p^2$, or $(pq-1)/p^2$. Our main tool is
 an observation about the sums of continued fraction coefficients: if $p>q>0$ are
relatively prime integers, then the sum of the coefficients in the continued
fraction expansion of $p/q$ is at most $p$. This provides a bound on the crosscap
number that we can in turn compare to the genus or crossing number.

The paper is organised as follows. After this introduction, we give a
brief overview of Teragaito's classification in Section 2 as well as some basic
results about continued fractions.  In Section 3, we apply these techniques to
prove  Theorem~\ref{propcross}. In Section 4, we outline the proof of 
Theorem~\ref{propgen}.

\section{Continued Fractions and Teragaito's classification}

In this section, we give an overview of Teragaito's~\cite{T} classification of 
the crosscap number of torus knots as well as some basic facts about
continued fractions that will prove useful in the sequel.
Recall that a positive rational number $r$ can be represented uniquely
by a simple continued fraction
$$r = a_0 + \frac{1}{a_1 + \frac{1}{\ddots + \frac{1}{a_{n-1}+\frac{1}{a_n}}}}$$
where each $a_i$, $1 \leq i \leq n$, is a positive integer and $a_n > 1$. We 
will write $r = {[} a_0, a_1, \ldots a_n {]}$.
If $\frac{p}{q} = {[} a_0, a_1, \ldots a_n {]}$, then
the function $N(p,q)$ defined by Bredon and Wood~\cite{BW} is given as follows.
We first sum the $a_i$ in order, beginning with $a_0$, and skipping the succeeding
$a_i$ whenever the partial sum becomes even. We then halve the total. For example,
$8/3 = {[} 2,1,2 {]}$, so $N(8,3) = (2 + 2)/2 = 2$. Since 
$34/49 = {[} 0,1,2,3,1,3 {]}$, we have $N(34,49) = (0 + 2 + 1 + 3)/2 = 3$.

For a $(p,q)$ torus knot $K$ with $p, q > 0$, we will say $K$ is odd (respectively,
even) if $pq$ is odd (respectively, even). For the statement of Teragaito's
theorem, we will assume $p > q$ if $K$ is odd and $p$ is even if $K$ is even.

\begin{thm}[Theorem 1 of \cite{T}]
Let $K$ be the non-trivial torus knot of type $(p,q)$, where $p,q > 0$.
\begin{enumerate}
\item If $K$ is even, then $c(K) = N(p,q)$.
\item If $K$ is odd, then $c(K) = \min \{ N(pq -1, p^2), N(pq+1,p^2) \}$
\end{enumerate}
\label{thmTer}
\end{thm}

In light of Theorem~\ref{thmTer}, the following observation about
continued fractions will be useful in bounding the crosscap numbers of torus knots.
We omit the straightforward proof by induction.
\begin{lem} Let $p > q > 0$ be relatively prime and let 
$p/q = {[}a_0, a_1, \ldots a_n {]}$ be a simple continued fraction.
Then $\sum_{i=0}^n a_i \leq p$.
\label{lemstandard}
\end{lem}

Finally, we will make use of a lemma that relates the continued fractions
of $(pq \pm 1)/p^2$ and $q/p$. For this lemma, assume $p > q > 1$.
\begin{lem}[Lemma 9 of \cite{T}]
If $q/p = {[} a_0, a_1, a_2, \ldots, a_n {]}$, then 
$$(pq-1)/p^2 = \left\{ \begin{array}{l}
{[} a_0, a_1, a_2, \ldots, a_{n-1}, a_n + 1, a_n - 1, \\
\mbox{ } a_{n-1}, a_{n-2}, \ldots, a_2, a_1 {]} \mbox{ if } n \mbox{ is odd, } \\
{[} a_0, a_1, a_2, \ldots, a_{n-1}, a_n - 1, a_n + 1, \\
\mbox{ } a_{n-1}, a_{n-2}, \ldots, a_2, a_1 {]} \mbox{ if } n \mbox{ is even, } \\
\end{array} \right.$$
and 
$$(pq+1)/p^2 = \left\{ \begin{array}{l}
{[} a_0, a_1, a_2, \ldots, a_{n-1}, a_n - 1, a_n + 1, \\
\mbox{ } a_{n-1}, a_{n-2}, \ldots, a_2, a_1 {]} \mbox{ if } n \mbox{ is odd, } \\
{[} a_0, a_1, a_2, \ldots, a_{n-1}, a_n + 1, a_n - 1, \\
\mbox{ } a_{n-1}, a_{n-2}, \ldots, a_2, a_1 {]} \mbox{ if } n \mbox{ is even. } \\
\end{array} \right.$$
\label{lem9}
\end{lem}

Note that, since $p>q$, $a_0 = 0$ in the lemma. Also, if $a_1 = 1$, then
${[} a_0, \ldots a_3,a_2,a_1 {]}$ should be replaced by 
${[} a_0, \ldots a_3,a_2 + 1{]}$.

\section{Proof of Theorem~\ref{propcross}}

\setcounter{thm}{1}

In this section, we will prove
\begin{thm} For a torus knot $K$, $c(K) \leq \lfloor (n(K) + 16)/12 \rfloor$.
\end{thm}

\Pf
By definition~\cite{C}, the crosscap number of the unknot is zero and the
theorem holds in this case. 
So, let $K$ be a $(p,q)$ torus knot where $p>q>1$ are relatively prime integers. 
The proof breaks into three cases according to whether $q$ is even, 
$p$ is even, or both are odd.

\subsection{$q$ is even}

Let us assume $q$ is even. By the Euclidean algorithm, there are unique positive
integers $m$ and
$k$ with $p = qm - k$ and $k < q$. Since $p > q$, we have $m > 1$. The crosscap
number of $K$,  
$N(q,p)$, is determined by the continued fraction
$$ q/p = q/(qm-k) = 0 + \frac{1}{(m-1) + \frac{1}{\frac{q}{q-k}}}$$
In calculating $N(p,q)$ we would skip $m-1$ (since the first
partial sum $a_0 = 0$ is even) and add certain of the coefficients
in the continued fraction of $q/(q-k)$. 
By Lemma~\ref{lemstandard}, the sum of all the coefficients in 
$q/(q-k)$ is bounded by $q$. For $N(p,q)$, the sum is halved, so 
we have $c(K) = N(p,q) \leq q/2$. 

Since $q < p$, the crossing number of $K$ is $n(K) = p(q-1) = (qm-k)(q-1)$.
In order to prove the theorem in this case, it's enough to show that 
\begin{equation}
\frac{q}{2} \leq \frac{(qm-k)(q-1) + 16}{12}. \label{eqnqeven}
\end{equation}

Our strategy is to use induction on $m$ and $k$. In the $m$ inductive step, $m$ 
will increase by one, while for the $k$ induction, we'll decrement by one at each
step.  We have already mentioned that $m > 1$ and, since $k
< q$, we begin our induction with the case
$m=2$ and $k = q-1$.
Then Equation~\ref{eqnqeven} becomes $q^2 -6q + 15 \geq 0$ which is true
for all integers $q$. So the theorem is proved in this case. 
If $m$ is increased by $1$, then the right hand side of Equation~\ref{eqnqeven} is
increased by $q(q-1)/12$ while the left hand side is unchanged. Since $q > 1$, the
equation will still hold if $m$ is increased by $1$. For the $k$ induction,
if $k$ is decreased by $1$, the right
hand side of Equation~\ref{eqnqeven} is increased by $(q-1)/12$ while the left
hand side is unchanged. By induction, Equation~\ref{eqnqeven} holds for all $m >1$
and all $0 \leq k < q$. This proves the theorem in the case $q$ is even.

\subsection{$p$ is even}

Let us assume $p$ is even. In this case $c(K) = N(p,q)$.
We can write $p = 2qm - k$ for some positive integers
$m$ and $k$ with $k < 2q$. If $k = q$, then $q \mid p$ contradicting our
assumption that $p$ and $q$ are relatively prime. We have two subcases depending
on whether $k < q$ or $k > q$.

Suppose $k < q$. The continued fraction for $N(p,q)$ is
$$p/q = (2qm-k)/q = 2m-1 + \frac{1}{\frac{q}{q-k}}$$ so that
$c(K) = N(p,q) \leq (2m-1 + q)/2$ (using Lemma~\ref{lemstandard}). On the other
hand, the crossing number is $n(K) = p(q-1) = (2qm-k)(q-1)$, so the theorem
can be proved by verifying
\begin{equation}
\frac{2m-1 + q}{2} \leq \frac{(2qm-k)(q-1) + 16}{12}.
\label{eqnpeven}
\end{equation}

Again, we will use induction on $m$ and $k$. If $m = 1$ and $k = q-1$, the
inequality becomes $(q-3)^2 \geq 0$. Note that, as $p$ is even, $q$ is odd. We
were already assuming $q>1$, so we must have $q \geq 3$. Thus,
Equation~\ref{eqnpeven} holds in the case $m = 1$, $k = q-1$. If $m$ is increased 
by $1$, the left hand side of the inequality is increased by $1$ while the
right is increased by $2q(q-1)/12$. Since $q \geq 3$, we have $2q(q-1)/12 \geq
1$ and the inductive step for $m$ is proved. If $k$ is decreased by $1$,
the left hand side is unchanged while the right hand side increases by
$(q-1)/12$. Thus, Equation~\ref{eqnpeven} holds for all $m \geq 1$ and
$k < q$ and the theorem is proved in the case where $k < q$.

If $k > q$, we have 
$$p/q = (2qm-k)/q = 2m-2 + \frac{1}{\frac{q}{2q-k}}$$
so that $c(K)  = N(p,q) \leq m-1 + q/2$. In this case we must verify the inequality
\begin{equation}
m-1 + \frac{q}{2} \leq \frac{(2qm-k)(q-1) + 16}{12} 
\label{eqn2peven}
\end{equation}
If $m=1$, then $p = 2qm-k = 2q-k < q$ which contradicts our assumption
that $p>q$. Therefore, the base step for the induction is $m=2$ and 
$k = 2q-1$. With these values, Equation~\ref{eqn2peven} becomes
$2q^2-7q+3 \geq 0$ and this inequality holds for all $q \geq 3$. The induction
for $m$ and $k$ is similar to the previous subcase and, thus, the theorem
is  proved for all even $p$.

\subsection{$pq$ odd}

Suppose both $p$ and $q$ are odd, and
let $p = qm -k$ where $m$ and $k$ are positive integers with $k < q$. Since
$p>q$, we have $m > 1$.
In this case $c(K)$ is determined by the continued fractions of $(pq \pm 1)/p^2$
and, by Lemma~\ref{lem9}, these are related to $q/p$.
Now, $$q/p = q/(qm-k) = 0 + \frac{1}{m-1 + \frac{1}{\frac{q}{q-k}}}.$$
So, if we write $q/p = {[} a_0, a_1, \ldots a_n {]}$ as in
Lemma~\ref{lem9}, then $a_0 = 0$, $a_1 = m-1$, and $q/(q-k) = {[} a_2, a_3,
\ldots,  a_n {]}$. Moreover, by Lemma~\ref{lemstandard}, 
$\sum_{i=2}^n a_i \leq q$.

For odd knots, $c(K) = \min \{ N(pq-1,p^2), N(pq+1,p^2) \}$. By Lemmma~\ref{lem9},
the continued fractions for $(pq \pm 1)/p^2$ both begin ${[} 0, m-1, a_2, \ldots
{]}$. So in calculating $N(pq \pm 1,p^2)$ we omit the $m-1$ coefficient and
the summation effectively begins with $a_2$. Moreover, by the lemma, 
whichever of $(pq+1)/p^2$ and $(pq-1)/p^2$ is used (and whether or not
$n$ is even),
the sum of the coefficients beginning with $a_2$ is $2 \sum_{i=2}^n a_i + a_1$.
Thus, $c(K) \leq \sum_{i=2}^n a_i + a_1/2 \leq q + (m-1)/2$. Since $c(K)$ is 
an integer, we have $c(K) \leq \lfloor q + (m-1)/2 \rfloor$.

The crossing number in this case is again $n(K) = p(q-1)$; so we can prove the
theorem by verifying the inequality:
\begin{equation} 
\lfloor q + \frac{m-1}{2} \rfloor \leq \frac{(qm-k)(q-1) + 16}{12}.
\label{eqnodd}
\end{equation}
In fact, this inequality does not hold for all choices of $q$, $m$, and $k$. 
Let us begin by delineating the cases where it does hold. 

Since $p$ and $q$ are both odd, we will need to carry out two induction arguments,
one for the case where $m$ is odd and $k$ is even and one with the opposite
parities. The base case for $m$ even is $m = 2$, $k = q-2$. In this case
Equation~\ref{eqnodd} becomes $q^2-11q + 14 \geq 0$ which is valid for 
all $q > 9$. The base case for $m$ odd is $m = 3$, $k = q-1$. Here the inequality
becomes $2q^2-13q + 3 \geq 0$ which is valid for all $q \geq 7$. Since $q$ is
odd, the smallest $q$ for which both cases apply is $q=11$. Let us
show the induction for $q \geq 11$ and then examine smaller values of $q$
individually.

\noindent%
\underline{Case 1} $q \geq 11$.
We've established that the base cases both hold if $q \geq 11$. We're left to
verify the inductive steps. Note that in both inductions $m$ will be increased
by $2$ at each step and $k$ will be decreased by $2$. If $m$ is increased
by $2$, the left hand side of Equation~\ref{eqnodd} is increased by $1$
while the right hand side increases by $q(q-1)/6$. 
Since $q \geq 11$, the
$m$ inductive step preserves the inequality. 
(Indeed, this 
inductive step will be valid so long as $q \geq 3$.) 
If $k$ is decreased by $2$,
the right hand side is increased by $(q-1)/6$ and the left hand side is
unchanged. Thus, the theorem is proved in the case $q \geq 11$.

\noindent%
\underline{Case 2} $q = 9$. Since the $m=2$ base case is problematic,
let's instead begin the even $m$ induction with $m = 4$ and
$k = q-2$. Then Equation~\ref{eqnodd} becomes $3q^2 - 13q + 2 \geq 0$ which is true
for all $q \geq 5$. As noted above, the $m=3$ base case is valid, and the inductive
arguments also go through when $q=9$. So, in order to complete this case,
we must address the knots that have $m=2$. That is, we must verify the
theorem for the knots $(17,9)$, $(13,9)$, and $(11,9)$. (Note that 
$15$ and $9$ are not relatively prime.)
The crossing numbers $n(K)$ of these knots are, respectively, $136$, 
$104$, and $88$. The crosscap numbers $c(K)$ are $5$, $4$, and $5$.
Thus, the theorem holds in these cases as well and is proved for the
case $q = 9$.

\noindent%
\underline{Case 3} $q = 7$. Our arguments above show that
the theorem holds for $q =7$ provided $m \geq 3$. Again, we
can verify the knots with $m = 2$, $(13,7)$, $(11,7)$, and $(9,7)$, directly.
The crossing numbers for these knots are $78$, $66$, and $54$ while the 
crosscap numbers are $4$, $3$, and $4$. So these knots also satisfy the theorem.

\noindent%
\underline{Case 4} $q = 5$. When $q=5$, the $m=3$ base case is no longer
valid. However, 
Equation~\ref{eqnodd} is satisfied if we take $q = 5$, $m = 5$ and $k = q-1 = 4$.
Thus, induction arguments will take care of all cases where $m \geq 4$. We
are left to investigate $m = 2$ and $m=3$. That is, we are left with 
the knots $(13,5)$, $(11,5)$, $(9,5)$, and $(7,5)$. The crossing
numbers are $52$, $44$, $36$, and $28$ while the crosscap numbers are 
all $3$. This completes the argument in the case $q = 5$.

\noindent%
\underline{Case 5} $q = 3$. If $q = 3$, we can explicitly calculate the continued
fraction $q/p$. Since $p$ is odd and relatively prime to $3$, $p$ is of the form
$6m+1$ or $6m-1$.

If $p = 6m+1$, then $q/p = {[} 0, 2m, 3 {]}$. As Teragaito~\cite{T} shows,
since $3x \equiv -1 \bmod p$ has the even solution $x = 2m$, the crosscap
number is $c(K) = N(pq-1,p^2)$. By Lemma~\ref{lem9},
$(pq-1)/p^2 = {[}0, 2m,2,4,2m {]}$. Therefore, $c(K) = (0 + 2 + 2m)/2 = m+1$.
On the other hand, the crossing number is $n(K) = p(q-1) = 2(6m+1)$.
Thus, 
$\lfloor (n(K) + 16)/12 \rfloor = 
\lfloor (12m+ 18)/12 \rfloor = m+1 = c(K)$ 
and the theorem holds in this case.

Finally, suppose $p = 6m-1$ and $q = 3$. Then,
$q/p = {[}0, 2m-1, 1,2 {]}$. Since $3x \equiv -1 \bmod p$ has the odd solution
$x = 4m-1$, the crosscap number is $N(pq+1,p^2)$ (see~\cite{T}). By
Lemma~\ref{lem9}, $(pq+1)/p^2 = {[}0,2m-1,1,1,3,1,2m-1{]}$ if $m \neq 1$
and ${[}0,1,1,1,3,2{]}$ when $m = 1$. Thus, $c(K) = (1+1+1+2m-1)/2 = m+1$. Now,
the crossing number is $p(q-1) = 2(6m-1)$. Thus,
$\lfloor (n(K) + 16)/12 \rfloor = \lfloor (12m +14)/12 \rfloor = m+1 = c(K)$
and the theorem holds in this case as well. 

Thus, we have proved Theorem~\ref{propcross} when $pq$ is odd. This completes 
the proof of the theorem. \qed

\section{Proof of Theorem~\ref{propgen}}

\setcounter{thm}{0}

In this section we prove 
\begin{thm} 
 For a torus knot $K$, $c(K) \leq \lfloor (g(K)+9)/6 \rfloor$.
\end{thm}

\Pf
The argument is very similar to that used in proving Theorem~\ref{propcross}.
The theorem holds for the trivial knot, so we will assume $p > q > 1$ are
relatively prime and $K$ is the $(p,q)$ torus knot. We have three cases depending
on whether $q$ is even, $p$ is even, or $pq$ is odd.

Let's assume $q$ is even. Then as in the previous section,
we have $c(K) \leq q/2$. Writing $p = qm - k$ with $m > 1$ and $k <q$ positive
and using $g(K) =  (p-1)(q-1)/2$, we can prove the theorem by verifying
\begin{equation}
\frac{q}{2} \leq \frac{(qm-k-1)(q-1)/2 +9}{6}
\label{eqngq}
\end{equation}
If $m=2$ and $k = q-1$, the inequality becomes $q^2 -7q + 18 \geq 0$ which is 
true for all integers $q$. If $m$ is increased by one or $k$ is decreased by
one, the right hand side increases while the left hand side is unchanged. 
Thus the inequality holds for all $m > 2$ and all positive $k < q$. This proves
the theorem in the case $q$ is even.

Next, assume $p$ is even and write $p = 2qm - k$ with $m$ and $k <2q $
positive. We have two subcases: $k < q$ and $k > q$. Suppose $k < q$. Then, as in
the previous section, $c(K) \leq (2m - 1 + q)/2$ and we can prove the theorem
by verifying
\begin{equation}
\frac{2m-1 + q}{2} \leq \frac{(2qm-k-1)(q-1)/2 +9}{6}.
\label{eqngp1}
\end{equation}
If $m=1$ and $k = q-1$, the inequality becomes $(q-3)(q-4) \geq 0$ which is true
for all integers $q$. Increasing $m$ by one will increase the left hand side
by one and  the right hand side by $q(q-1)/6$. So, this inductive step will go
through for all $q \geq 3$. If $k$ is decreased by one, the left hand side is
unchanged and the right hand side increases so the $k$ induction step will also
preserve the inequality. So the theorem is proved when $k < q$. 

%Now, $p$ is even, so $q$ must be odd and $q > 1$. So we are left to consider the 
%case $q = 3$. If $m=2$ and $k = q-1$, Equation~\ref{eqngp1} becomes 
%$q(q-3) \geq 0$. So, using this as a base case, we can apply induction to prove
%the theorem for all $m \geq 2$ when $q = 3$. Finally, if $q=3$, $m = 1$,
%and $k < q$, we are talking about the knot $(4,3)$. This knot has genus
%$g(K) = 3$ and crosscap number $c(K) = 2$ so the theorem is proved for all
%$q$ when $k < q$.

Now suppose $k > q$ and $p = 2qm-k$ is even. As in the previous
section, we have $c(K) \leq m-1 + q/2$ so it 
will be enough to show
\begin{equation}
m-1 + \frac{q}{2} \leq \frac{(2qm-k-1)(q-1)/2 +9}{6}.
\label{eqngp2}
\end{equation}
If $m = 1$, since $k < q$, then $p = 2qm-k$ will be less than $q$ contradicting an
earlier assumption. So $m \geq 2$. Substituting $m=2$ and $k = 2q-1$ into
Equation~\ref{eqngp2} results in the inequality $(q-1)(q-3) \geq 0$ which is true
since $q \geq 3$. Again, the $m$ and $k$ inductive steps preserve the
inequality and the theorem is proved in the case $k <q$ as well.

Finally, we have the case where $pq$ is odd. As in the previous section,
we can show $c(K) \leq \lfloor q + (m-1)/2 \rfloor$. It's enough
to verify 
\begin{equation}
\lfloor q + \frac{m-1}{2} \rfloor \leq \frac{(qm-k-1)(q-1)/2 +9}{6}.
\label{eqngodd}
\end{equation}

If $m = 2$ and $k = q-2$ we have $q^2 - 12q + 17 \geq 0$ which is valid for
$q \geq 11$. If $m = 3$ and $k = q-1$, we have $2q^2 - 14q + 6 \geq 0$ which is
valid for $q \geq 7$. So we will look at the induction when $q \geq 11$ and
then take smaller values of $q$ in turn.

\noindent%
\underline{Case 1} $q \geq 11$. If $q \geq 11$, both base cases hold, so it's
enough to check that the inductive steps preserve the inequality. If $m$ is
increased by $2$, the left side of Equation~\ref{eqngodd} goes up by $1$
while the right hand side increases by $q(q-1)/6$. Thus, the $m$ induction
will work as long as $q \geq 3$. If $k$ is decremented by $2$, the left hand
side of the inequality is unchanged while the right increases by $(q-1)/6$. 
Therefore, the theorem is proved when $q \geq 11$.

\noindent%
\underline{Case 2} $q = 9$. If $m = 4$ and $k = q-2$, Equation~\ref{eqngodd}
becomes $3q^2-14q+5 \geq 0$ which is valid for $q \geq 5$. So, induction
arguments show the theorem holds for $m \geq 3$. For $m = 2$, we have
the knots $(17,9)$, $(13,9)$, and $(11,9)$ of genus $72$, $48$, and $40$
respectively. As these knots have crosscap number $5$, $4$, and $5$, the
theorem holds for these knots and, therefore, for all knots with $q = 9$.

\noindent%
\underline{Case 3} $q = 7$. As above, induction will take care of all cases
where $m \geq 3$. For $m = 2$, we have the knots $(13,7)$, $(11,7)$, and
$(9,7)$ of genus $36$, $30$, and $24$ respectively. Since the crosscap numbers
are $4$, $3$, and $4$, we have verified the theorem in the case $q = 7$.

\noindent%
\underline{Case 4} $q = 5$. If $m = 5$ and $k = q-1$, Equation~\ref{eqngodd}
becomes $4q^2 -16q -6 \geq 0$ which is valid when $q = 5$. So, induction
takes care of the cases where $m \geq 4$ and we're left with the knots
$(13,5)$, $(11,5)$, $(9,5)$, and $(7,5)$ of genus $24$, $20$, $16$, and $12$
respectively. These all have crosscap number $3$, so the theorem holds
in this case as well.

\noindent%
\underline{Case 5} $q = 3$. If $p = 6m+1$, then $c(K) = m+1$,
$g(K) = 6m$, and the theorem holds. If $p = 6m-1$, then $c(K) = m+1$,
$g(K) = 6m-2$, and the theorem holds. 

This completes the proof of Theorem~\ref{propgen}. \qed

\end{document}